\documentclass[11pt]{article}

\date{11th February 2021}

\title{\vskip-1.0em\sc Unavoidable subprojections in union-closed set systems of infinite breadth}
\author{\sc Yemon Choi, Mahya Ghandehari, Hung Le Pham}

\RequirePackage{amsmath,amssymb, graphicx, mathtools}

\usepackage[heavycircles]{stmaryrd}

\usepackage{tikz}
\usepackage{tkz-euclide}

\RequirePackage{amsthm}   
\newcounter{pulse}[section]
\numberwithin{pulse}{section}  

\newcommand{\thf}{\sc} 

\theoremstyle{plain}
\newtheorem{thm}[pulse]{\thf Theorem}

\newtheorem{prop}[pulse]{\thf Proposition}
\newtheorem{lem}[pulse]{\thf Lemma}

\theoremstyle{definition}
\newtheorem{dfn}[pulse]{\thf Definition}
\newtheorem*{notn}{\thf Notation}
\newtheorem{eg}[pulse]{\thf Example}
\theoremstyle{remark}
\newtheorem{rem}[pulse]{\thf Remark}

\newenvironment{romnum}{%
\begin{enumerate}

}{\end{enumerate}\ignorespacesafterend}

\newcommand{\gm}{\gamma}

\newcommand{\Gm}{\Gamma}
\newcommand{\Om}{\Omega}

\newcommand{\defeq}{:=}
\newcommand{\dt}[1]{\textit{#1}\/}  
\newcommand{\st}{\mathbin{\colon}}

\newcommand{\join}{\operatorname{join}}
\newcommand{\meet}{\operatorname{meet}}

\newcommand{\Inc}{\operatorname{Inc}}

\newcommand{\FIN}{\operatorname{\mathcal{P}}\nolimits^{\rm fin}} 
\newcommand{\COFIN}{\operatorname{\mathcal{P}}\nolimits^{\rm cofin}} 

\newcommand{\ssfont}[1]{{\mathsf{#1}}}
\newcommand{\ssa}{\ssfont{a}}
\newcommand{\ssb}{\ssfont{b}}

\newcommand{\ssp}{\ssfont{p}}

\newcommand{\ssw}{\ssfont{w}}
\newcommand{\ssx}{\ssfont{x}}
\newcommand{\ssy}{\ssfont{y}}
\newcommand{\ssz}{\ssfont{z}}

\newcommand{\cE}{{\mathcal E}}
\newcommand{\cF}{{\mathcal F}}
\newcommand{\cG}{{\mathcal G}}
\newcommand{\cH}{{\mathcal H}}
\newcommand{\cI}{{\mathcal I}}
\newcommand{\cJ}{{\mathcal J}}
\newcommand{\cL}{{\mathcal L}}
\newcommand{\cP}{{\mathcal P}}

\newcommand{\cS}{{\mathcal S}}
\newcommand{\cT}{{\mathcal T}}
\newcommand{\cU}{{\mathcal U}}
\newcommand{\cV}{{\mathcal V}}

\newcommand{\Nat}{{\mathbb N}}

\newcommand{\TMAX}{{\cT}_{\rm max}}
\newcommand{\TMIN}{{\cT}_{\rm min}}
\newcommand{\TORT}{{\cT}_{\rm ort}}

\newcommand{\MAX}{{\rm MAX}}
\newcommand{\MIN}{{\rm MIN}}
\newcommand{\ORT}{{\rm ORT}}

\newcommand{\br}{\operatorname{br}} 

\newcommand{\subproj}{\mathop{\leq_{\sf p}}} 

\makeatletter
\newcommand{\optfs}{\@ifnextchar.{}{.\@}}
\makeatother
\newcommand{\ucib}{u.c.i.b\optfs}

\newcommand{\sscap}{\wedge} 
\newcommand{\sscup}{\vee} 

\newcommand{\bigsscap}{\bigwedge} 
\newcommand{\bigsscup}{\bigvee} 

\newcommand{\setcap}{\cap} 
\newcommand{\setcup}{\cup} 
\newcommand{\bigsetcap}{\bigcap} 
\newcommand{\bigsetcup}{\bigcup} 

\newcommand{\dotcup}{\mathbin{\dot{\cup}}}

\usepackage[a4paper,left=30mm,right=30mm,top=30mm,bottom=30mm]{geometry}

\usepackage[colorlinks=true,linkcolor=blue,citecolor=red]{hyperref}
\usepackage{sectsty}
\allsectionsfont{\sffamily}

\begin{document}

\maketitle

\begin{abstract}
We consider union-closed set systems with infinite breadth, focusing on three particular configurations ${\mathcal T}_{\rm max}({\mathcal E})$, ${\mathcal T}_{\rm min}({\mathcal E})$ and ${\mathcal T}_{\rm ort}({\mathcal E})$. We show that these three configurations are not isolated examples; in any given union-closed set system of infinite breadth, at least one of these three configurations will occur as a subprojection.
This characterizes those union-closed set systems which have infinite breadth,
and is the first general structural result for such set systems.
\\
 

\end{abstract}

\setcounter{page}{0}

\bigskip\hrule
\paragraph{Notification:}
{\sf this update of the arXiv version fixes two minor errors in Propositions \ref{p:reformulate TMIN} and~\ref{p:reformulate TORT}, which were spotted after the paper was accepted and were fixed during the ``page proofs'' stage.
In both cases the errors were in the implications (i)$\implies$(ii), and did not affect the validity of the proof of the main theorem (Theorem \ref{t:headline result}), which only requires the implications (ii)$\implies$(i).

Nevertheless, for the sake of accuracy we have decided to upload a corrected version to arXiv. We have also updated the bibliography: the paper \cite{CGP1} has now been published.

This update of the arXiv version does not include the reformatting and copy-editing done by the journal, and does not supersede the official ``version of record'', which has been published online as

\begin{quote}
{\sc Y. Choi, M. Ghandehari, H. L. Pham}, 
Unavoidable subprojections in union-closed set systems of infinite breadth.
European J. Combin.~94 (2021), article 103311.
DOI \href{https://doi.org/10.1016/j.ejc.2021.103311}{10.1016/j.ejc.2021.103311}
\end{quote}

\noindent
Note that the version of record includes the corrections mentioned above.
}
\bigskip\hrule
\newpage

\tableofcontents

\bigskip\hrule\bigskip

\vfill\eject

\begin{section}{Introduction}
\label{s:intro}

\begin{subsection}{Initial definitions}
\label{ss:initial def}
For a given set $\Om$, a \dt{set system on $\Om$} is a subset of the powerset $\cP(\Om)$.
Usually we study set systems with extra structure; this paper concerns those set systems 
which satisfy the following additional property.

\begin{dfn}\label{d:union-closed}
Let $\cS$ be a set system on $\Om$. We say that $\cS$ is \dt{union-closed} if $\ssa,\ssb\in \cS\implies \ssa\setcup \ssb\in \cS$.
\end{dfn}

For such an $\cS$, the binary operation $(\ssa,\ssb)\mapsto \ssa\setcup\ssb$ makes $\cS$ into a commutative semigroup in which every element is idempotent, a so-called \dt{(join) semilattice.}
Hence, in studying union-closed set systems, we may import concepts from lattice theory.\footnote{For instance, it is well known that Frankl's conjecture on \emph{finite} union-closed set systems has a lattice-theoretic reformulation, which has then been solved for certain natural classes of lattices.}
One such lattice-theoretic concept is the notion of \dt{breadth}.

\begin{dfn}
\label{d:breadth}
Let $\cS\subseteq\cP(\Om)$ be a union-closed set system. The \dt{breadth of $\cS$}, denoted by $\br(\cS)$, is the smallest $n\in\Nat$ with the following property: whenever $\ssx_1,\dots, \ssx_{n+1}\in \cS$, there exists $j\in\{1,\dots,n+1\}$ such that
\[
\bigsetcup_{i=1}^{n+1} \ssx_i = \bigsetcup_{i\neq j} \ssx_i
\]
If no such $n$ exists, we say that $\cS$ has \dt{infinite breadth}
and we put $\br(\cS)=\infty$.
\end{dfn}

It is easily checked from this definition that if $\Om$ is a non-empty \emph{finite set}, then $\br(\cP(\Om))=|\Om|$; the lower bound on breadth is demonstrated by considering singleton subsets of $\Om$, and the upper bound follows from the pigeonhole principle. The same reasoning shows that if $\Om$ is infinite then $\cP(\Om)$ has infinite breadth.

To date, there has been no attempt to develop any kind of structure theory for union-closed set systems with infinite breadth.
Informally, if $\cS$ is a union-closed set system with infinite breadth, then for each $n\in\Nat$ there exists  a finite subset $\cF_n\subset \cS$ which is at least as complex as the powerset of $\{1,\dots,n\}$.
However, the condition of infinite breadth gives us no {\it a priori} information or constraint on how such $\cF_n$ fit together.
 The starting point for the present paper is the observation that there are three particular set systems, union-closed with infinite breadth, where the subsets $\cF_n$ are positioned relative to each other in a very structured way.

\begin{dfn}[Three key configurations with infinite breadth]
\label{d:3kings}
Let $\Gamma$ be a countably infinite set, and let $\cE$ denote a partition of $\Gamma$ into disjoint subsets $(E_n)_{n\geq 1}$ with $|E_n|=n+1$ for all $n$.
Let
 $E_{<n} = \bigsetcup_{1\leq j \leq n-1} E_j$ and
 $E_{>n} = \bigsetcup_{j \geq n+1} E_j$, with the convention that $E_{<1}=\emptyset$.

We define the following union-closed set systems on $\Gamma$ (see Figure \ref{fig1} for an illustration):
\begin{romnum}
\item 
$\TMAX(\cE)  \defeq \{ E_{<n} \cup \ssx \colon n\in\Nat, \ssx\subseteq E_n, \emptyset\neq \ssx \}$;
\item
$\TMIN(\cE) \defeq \{ E_{>n} \cup \ssx \colon n\in\Nat, \ssx\subseteq E_n, \emptyset\neq \ssx \}$;
\item
$\TORT(\cE) \defeq \{ E_{<n}\cup E_{>n} \cup \ssx \colon n\in\Nat, \ssx\subseteq E_n, \emptyset\neq \ssx \}$.
\end{romnum}
\end{dfn}

\begin{figure}[hpt]
\centering
\begin{tikzpicture}[scale=0.3]
	\tkzLabelPoint[right,scale=0.6](-3,0.4){$E_1$};
	\draw[black, pattern=north east lines, pattern color=black][-] (0,0)--(2,0) --(2,1)--(0,1)--cycle;
	\tkzLabelPoint[right,scale=0.6](-3,-0.6){$E_2$};
	\draw[black, pattern=north east lines, pattern color=black][-] (0,-1)--(3,-1) --(3,0)--(0,0)--cycle;
	\tkzLabelPoint[scale=0.4](1,-1.2){$\bullet$};
	\tkzLabelPoint[scale=0.4](1,-1.8){$\bullet$};
	\tkzLabelPoint[scale=0.4](1,-2.4){$\bullet$};
	\tkzLabelPoint[right,scale=0.6](-3,-3.6){$E_{n-2}$};
	\draw[black, pattern=north east lines, pattern color=black][-] (0,-4)--(6,-4) --(6,-3)--(0,-3)--cycle;
	\tkzLabelPoint[right,scale=0.6](-3,-4.6){$E_{n-1}$};
	\draw[black, pattern=north east lines, pattern color=black][-] (0,-5)--(7,-5) --(7,-4)--(0,-4)--cycle;
	\tkzLabelPoint[right,scale=0.6](-3,-5.6){$E_{n}$};
	\draw[black, pattern=north east lines, pattern color=black][-] (0,-6)--(1,-6) --(1,-5)--(0,-5)--cycle;
	\draw[black, pattern=north east lines, pattern color=black][-] (3,-6)--(5,-6) --(5,-5)--(3,-5)--cycle;
	\draw[black, pattern=north east lines, pattern color=black][-] (6,-6)--(8,-6) --(8,-5)--(6,-5)--cycle;
	\draw[-] (0,-6)--(8,-6) --(8,-5)--(0,-5)--cycle;
	\tkzLabelPoint[right,scale=0.6](-3,-6.6){$E_{n+1}$};
	\draw[-] (0,-7)--(9,-7) --(9,-6)--(0,-6)--cycle;
	\tkzLabelPoint[right,scale=0.6](-3,-7.6){$E_{n+2}$};
	\draw[-] (0,-8)--(10,-8) --(10,-7)--(0,-7)--cycle;
	\tkzLabelPoint[scale=0.4](1,-8.2){$\bullet$};
	\tkzLabelPoint[scale=0.4](1,-8.8){$\bullet$};
	\tkzLabelPoint[scale=0.4](1,-9.4){$\bullet$};
	\tkzLabelPoint[scale=1.1](1,-9.6){$\TMAX(\cE)$}
	\begin{scope}[shift={(16,0)}]
	\tkzLabelPoint[right,scale=0.6](-3,0.4){$E_1$};
	\draw[-] (0,0)--(2,0) --(2,1)--(0,1)--cycle;
	\tkzLabelPoint[right,scale=0.6](-3,-0.6){$E_2$};
	\draw[-] (0,-1)--(3,-1) --(3,0)--(0,0)--cycle;
	\tkzLabelPoint[scale=0.4](1,-1.2){$\bullet$};
	\tkzLabelPoint[scale=0.4](1,-1.8){$\bullet$};
	\tkzLabelPoint[scale=0.4](1,-2.4){$\bullet$};
	\tkzLabelPoint[right,scale=0.6](-3,-3.6){$E_{n-2}$};
	\draw[-] (0,-4)--(6,-4) --(6,-3)--(0,-3)--cycle;
	\tkzLabelPoint[right,scale=0.6](-3,-4.6){$E_{n-1}$};
	\draw[-] (0,-5)--(7,-5) --(7,-4)--(0,-4)--cycle;
	\tkzLabelPoint[right,scale=0.6](-3,-5.6){$E_{n}$};
	\draw[black, pattern=north east lines, pattern color=black][-] (0,-6)--(1,-6) --(1,-5)--(0,-5)--cycle;
	\draw[black, pattern=north east lines, pattern color=black][-] (3,-6)--(4,-6) --(4,-5)--(3,-5)--cycle;
	\draw[black, pattern=north east lines, pattern color=black][-] (7,-6)--(8,-6) --(8,-5)--(7,-5)--cycle;
	\draw[-] (0,-6)--(8,-6) --(8,-5)--(0,-5)--cycle;
	\tkzLabelPoint[right,scale=0.6](-3,-6.6){$E_{n+1}$};
	\draw[black, pattern=north east lines, pattern color=black][-] (0,-7)--(9,-7) --(9,-6)--(0,-6)--cycle;
	\tkzLabelPoint[right,scale=0.6](-3,-7.6){$E_{n+2}$};
	\draw[black, pattern=north east lines, pattern color=black][-] (0,-8)--(10,-8) --(10,-7)--(0,-7)--cycle;
	\tkzLabelPoint[scale=0.4](1,-8.2){$\bullet$};
	\tkzLabelPoint[scale=0.4](1,-8.8){$\bullet$};
	\tkzLabelPoint[scale=0.4](1,-9.4){$\bullet$};
	\tkzLabelPoint[scale=1.1](1,-9.6){$\TMIN(\cE)$}	
	\end{scope}
	\begin{scope}[shift={(32,0)}]
	\tkzLabelPoint[right,scale=0.6](-3,0.4){$E_1$};
	\draw[black, pattern=north east lines, pattern color=black][-] (0,0)--(2,0) --(2,1)--(0,1)--cycle;
	\tkzLabelPoint[right,scale=0.6](-3,-0.6){$E_2$};
	\draw[black, pattern=north east lines, pattern color=black][-] (0,-1)--(3,-1) --(3,0)--(0,0)--cycle;
	\tkzLabelPoint[scale=0.4](1,-1.2){$\bullet$};
	\tkzLabelPoint[scale=0.4](1,-1.8){$\bullet$};
	\tkzLabelPoint[scale=0.4](1,-2.4){$\bullet$};
	\tkzLabelPoint[right,scale=0.6](-3,-3.6){$E_{n-2}$};
	\draw[black, pattern=north east lines, pattern color=black][-] (0,-4)--(6,-4) --(6,-3)--(0,-3)--cycle;
	\tkzLabelPoint[right,scale=0.6](-3,-4.6){$E_{n-1}$};
	\draw[black, pattern=north east lines, pattern color=black][-] (0,-5)--(7,-5) --(7,-4)--(0,-4)--cycle;
	\tkzLabelPoint[right,scale=0.6](-3,-5.6){$E_{n}$};
	\draw[black, pattern=north east lines, pattern color=black][-] (1,-6)--(3,-6) --(3,-5)--(1,-5)--cycle;
	\draw[black, pattern=north east lines, pattern color=black][-] (5,-6)--(6,-6) --(6,-5)--(5,-5)--cycle;
	\draw[-] (0,-6)--(8,-6) --(8,-5)--(0,-5)--cycle;
	\tkzLabelPoint[right,scale=0.6](-3,-6.6){$E_{n+1}$};
	\draw[black, pattern=north east lines, pattern color=black][-] (0,-7)--(9,-7) --(9,-6)--(0,-6)--cycle;
	\tkzLabelPoint[right,scale=0.6](-3,-7.6){$E_{n+2}$};
	\draw[black, pattern=north east lines, pattern color=black][-] (0,-8)--(10,-8) --(10,-7)--(0,-7)--cycle;
	\tkzLabelPoint[scale=0.4](1,-8.2){$\bullet$};
	\tkzLabelPoint[scale=0.4](1,-8.8){$\bullet$};
	\tkzLabelPoint[scale=0.4](1,-9.4){$\bullet$};
	\tkzLabelPoint[scale=1.1](1,-9.6){$\TORT(\cE)$}	
	\end{scope}
\end{tikzpicture}
\caption{Typical members of $\TMAX$, $\TMIN$, and $\TORT$ at level $n$.}
\label{fig1}
\end{figure}

Note that $\TMAX(\cE)$, $\TMIN(\cE)$ and $\TORT(\cE)$ all have infinite breadth.
For, if we let $\cT_*$ denote one of these three set systems, then for each $n\in\Nat$ we may enumerate $E_n$ as $\{\gm_1,\dots, \gm_{n+1}\}$ and then choose $\ssx_1,\dots, \ssx_{n+1}\in\cT_*$ such that $\ssx_j\cap E_n= \{\gm_j\}$ for $1\leq j\leq n+1$, which shows that $\br(\cT_*)> n$.

Informally speaking, the main result of our paper implies that these special configurations $\TMAX(\cE)$, $\TMIN(\cE)$ and $\TORT(\cE)$  are not isolated examples,
but are unavoidable when considering union-closed set systems of infinite breadth.
\end{subsection}

\begin{subsection}{Statement of the main result}
\label{ss:state main}

To state our main result precisely, some further notation and terminology is needed.

\begin{dfn}[Subprojections]\label{d:subprojection}
Given $\cS\subseteq\cP(\Om)$ and a non-empty $\Gm\subseteq\Om$, the \dt{projection of $\cS$ on $\Gm$} is defined to be the set system
\[
 \cS\owedge\Gm \defeq \{ \ssx\setcap\Gm \colon \ssx\in \cS\}\subseteq\cP(\Gm).
\]
If $\cS'\subseteq\cS$, then a set system of the form $\cT=\cS'\owedge\Gm$ is called a \dt{subprojection of $\cS$} (on $\Gm$). We denote this by $\cT\subproj\cS$ (with the choice of $\Gm$ made clear where necessary).
\end{dfn}

\begin{prop}\label{p:breadth of subproj}
Let $\cS$ and $\cT$ be union-closed set systems on $\Om$ with $\cT\subproj\cS$. Then $\br(\cT)\leq\br(\cS)$.
\end{prop}
\begin{proof}
If $\cS$ has infinite breadth there is nothing to prove.
Suppose $\cS$ has finite breadth and let $n=\br(\cS)$. 
By assumption, there exist $\cS'\subseteq\cS$ and $\Gamma\subseteq\Om$ such that $\cT=\cS'\owedge\Gm$.
Given $\ssy_1,\dots, \ssy_{n+1}\in \cT$, there exist $\ssx_1,\dots, \ssx_{n+1}\in \cS'\subseteq \cS$ such that $\ssy_i=\ssx_i\setcap\Gamma$ for all~$i$. Since $\br(\cS)=n$, by re-indexing if necessary we may assume that $\bigsetcup_{i=1}^{n+1}\ssx_i = \bigsetcup_{i=1}^n \ssx_i$; then $\bigsetcup_{i=1}^{n+1} \ssy_i = \bigsetcup_{i=1}^n\ssy_i$. This shows that $\br(\cT)\leq n$.
\end{proof}

In particular, let $\cS$ be a union-closed set system and suppose that $\cT_*\subproj\cS$, where $\cT_*$ denotes one of the three configurations $\TMAX(\cE)$, $\TMIN(\cE)$ or $\TORT(\cE)$ from Definition~\ref{d:3kings}. Then $\cS$ has infinite breadth, since $\cT_*$ does. Our main result is that the converse holds.

\begin{thm}[Main result]\label{t:headline result}
Let $\cS\subseteq\cP(\Om)$ be a union-closed set system
with infinite breadth. Then
$\cS$ has a subprojection equal to either $\TMAX(\cE)$, $\TMIN(\cE)$, or $\TORT(\cE)$. More precisely, there exist:
\begin{itemize}
\item a sequence $(E_n)_{n\geq 1}$ of pairwise disjoint subsets of $\Om$, satisfying $|E_n|=n+1$ for all $n\in\Nat$, and
\item a union-closed set system $\cS'\subseteq \cS$,
\end{itemize}
such that when we put $\cE=(E_n)_{n\geq 1}$ and $\Gamma= \bigsetcup_{n\geq 1} E_n$, the subprojection $\cS'\owedge\Gamma$ is equal to either $\TMAX(\cE)$, $\TMIN(\cE)$ or $\TORT(\cE)$.
\end{thm}

To our knowledge, Theorem \ref{t:headline result} is the first general structural result for union-closed set systems of infinite breadth, and most of the paper is devoted to its proof.
Note that while a given $\cS$ could have more than one of these three configurations as a subprojection, a little work shows that each of these three cannot contain (a copy of) the other two; so the list of unavoidable configurations in our theorem is best possible.

\begin{eg}[Two examples as illustrations of Theorem~\ref{t:headline result}]
\label{eg:fin-and-cofin}
Given an infinite set $\Om$, let $\FIN(\Om)$ be the set of all finite subsets of $\Om$ and let $\COFIN(\Om)$ be the set of all cofinite subsets of $\Om$.
Both $\FIN(\Om)$ and $\COFIN(\Om)$ are union-closed, and both have infinite breadth.

Fix a countably infinite subset $\Gm\subseteq\Om$, and
let $\cE$ be a partition of $\Gm$ with the properties described in Definition~\ref{d:3kings}.
Then $\TMAX(\cE)\subset \FIN(\Gm)\subseteq\FIN(\Om)$, and thus $\TMAX(\cE)$ is a subprojection of $\FIN(\Om)$ on $\Gm$. Moreover, if we define two set systems on $\Om$~by
\[ \cS_{\min} \defeq \{ \ssx\setcup (\Om\setminus\Gm) \colon \ssx\in \TMIN(\cE)\} \subset \COFIN(\Om),\]
\[ \cS_{\rm ort} \defeq \{ \ssx\setcup (\Om\setminus\Gm) \colon \ssx\in \TORT(\cE)\} \subset \COFIN(\Om), \]
then $\TMIN(\cE)=\cS_{\min}\owedge\Gm$ and
$\TORT(\cE)=\cS_{\rm ort}\owedge\Gm$. Thus both $\TMIN(\cE)$ and $\TORT(\cE)$ are subprojections of $\COFIN(\Om)$ on $\Gm$.

On the other hand, since $\COFIN(\Om)\owedge\Gm=\COFIN(\Gm)$, every subprojection of $\COFIN(\Om)$ on $\Gm$ is contained in $\COFIN(\Gm)$. Since every $\ssx\in\TMAX(\cE)$ is finite, this shows that $\TMAX(\cE)$ is not a subprojection of $\COFIN(\Om)$.
By similar reasoning, which we leave to the reader, one can show that neither $\TMIN(\cE)$ nor $\TORT(\cE)$ arise as subprojections of $\FIN(\Om)$.
\end{eg}

\end{subsection}

\begin{subsection}{Some context and connections}\label{ss:wider context}

\subsubsection{Relation to existing results}
\label{sss:existing results}

\paragraph{Extending previous results of the authors.}
Let $\cS\subseteq\FIN(\Om)$ be union-closed with infinite breadth.
An earlier combinatorial result of the authors (see \cite[Proposition 4.10]{CGP1}), when translated into the terminology of the present paper, states
 that $\cS$ has a subprojection of the form $\TMAX(\cE)$.
This result becomes an easy consequence of Theorem~\ref{t:headline result} and the fact that, as observed at the end of Example \ref{eg:fin-and-cofin}, neither $\TMIN(\cE)$ nor $\TORT(\cE)$ can arise as subprojections of~$\cS$.
As we shall see, the proof of Theorem~\ref{t:headline result} requires substantially more work and several new ideas compared with the earlier special case in~\cite{CGP1}.

\paragraph{Analogies with structural results in group theory.}
One may regard Theorem \ref{t:headline result} as similar in spirit to results in group theory, which say that the free group on two generators is not merely an example of a non-amenable group, but is contained in every non-amenable group that satisfies some mild additional hypotheses. (For instance, every non-amenable \emph{linear group} has a free non-abelian subgroup; this follows from the \dt{Tits alternative.})

Moreover, subprojections as defined in Definition~\ref{d:subprojection} have an algebraic interpretation, analogous to the notion of \dt{subquotients} in group theory or representation theory. Subquotients are the building blocks for theories of composition series (of groups or representations); we shall not attempt to develop a corresponding theory in this paper, although this could be useful if one seeks to
obtain a more detailed structure theory for union-closed set systems of infinite breadth.

\paragraph{Unavoidable configurations in finite set systems.}
It will follow from Lemma~\ref{l:subprojection} below that although we defined subprojections as projections of subsystems, this is equivalent to taking subsystems of projections.\footnotemark\
Thus Theorem \ref{t:headline result} is similar in spirit to the main result of \cite{BB05}, which shows that one of three special patterns must occur as a subprojection in any sufficiently large finite set system.
\footnotetext{The authors of \cite{BB05} use the terminology ``trace'' where we have used the word ``projection''; we feel that our choice is slightly clearer for the purposes of the present paper.}

\paragraph{Ramsey theory.}
There is a Ramsey-theoretic flavour to the statement of Theorem~\ref{t:headline result}, but it does not seem to follow from the standard Ramsey-type theorems for (hyper)graphs. Nevertheless, the nested inductive arguments that will be used to prove Theorem \ref{t:headline result} are reminiscent of similar arguments in infinitary Ramsey theory, and it might be interesting to pursue further links in future work.

\subsubsection{Breadth in other settings}
One measure of complexity for general set systems is \dt{Vapnik--Chervonenkis dimension}, which has played a seminal role in various aspects of statistical learning theory.
It turns out that when $\emptyset\in\cS$ and $\cS$ is union-closed, its VC-dimension coincides with its breadth (a fuller explanation is provided in Remark \ref{r:breadth and VC}). Attention has mostly focused on examples with finite VC-dimension, since these systems are the ones for which uniform learning guarantees can be devised. Nevertheless, the techniques used to prove Theorem \ref{t:headline result} may lead to a better understanding of \emph{how} things go wrong in settings with infinite VC-dimension.

Breadth of set systems also makes an unexpected appearance in some recent work on model theory \cite{ADHMS_VC1,ADHMS_VC2}. It turns out that \emph{stable theories} --- such as the theory of algebraically closed fields, or the theory of a module over a given ring --- will naturally give rise to set systems with infinite breadth.
See \cite[Proposition 2.20]{ADHMS_VC1} and the surrounding remarks\footnotemark\ for further details.
\footnotetext{In the statement of \cite[Proposition 2.20]{ADHMS_VC1}, the assumption ``${\rm vc}(\Phi)>0$'' refers not to VC-dimension but to the related notion of \dt{VC-density}, which we shall not discuss here. The reader should also beware that in \cite[Section 2.4]{ADHMS_VC1}, breadth is defined using intersections rather than unions, which yields an equivalent but dual perspective to the one in this paper.}

We already remarked that every union-closed set system can be viewed as a semilattice. Conversely, every semilattice $S$ can be naturally represented as an intersection-closed set system (this appears to be folklore, but an explicit construction is sketched in \cite[Section 2.1]{CGP1}).
Hence, by taking complements, $S$ can be modelled by a union-closed set system~$\cS$; moreover, $S$ has infinite breadth in the sense of lattice theory if and only if $\cS$ has infinite breadth in the sense of Definition~\ref{d:breadth}.
While there has been work on semilattices with small breadth, such as \cite{Ditor,Gierz,Lawson71}, to our knowledge there has been no systematic study of semilattices with infinite breadth.
The tools developed in this paper to prove Theorem \ref{t:headline result} may be useful in further study of the structure of such semilattices.

\subsubsection{Applications to Banach algebras}
\label{sss:BA applications}
Given a semilattice $S$, one can use a representation of $S$ as a union-closed set system to construct submultiplicative weight functions on~$S$, which give rise to examples of Banach function algebras.
In recent work \cite{CGP1}, we obtained a combinatorial characterization of those weights for which the corresponding Banach algebras have the so-called AMNM property of \cite{BEJ_AMNM1}. Using \cite[Theorem 3.7]{CGP1} and Theorem~\ref{t:headline result} of the present paper, together with further combinatorial arguments, one can construct on \emph{every} semilattice of infinite breadth a weight whose corresponding Banach algebra fails to have the AMNM property;
this resolves a question which was raised implicitly by the first author in \cite{YC_wtsl-amnm} and pursued further in~\cite{CGP1}.
Details will appear in the article \cite{CGP3}.
\end{subsection}

\begin{subsection}{Overview of the proof of the main result}
\label{ss:overview of proof}
This subsection is intended to explain the strategy/intuition behind the technical arguments and definitions in the rest of the paper.
However, none of the material in this subsection is logically necessary for the proofs in the following sections. The reader who prefers to see precise details and proofs immediately, should skip to Section \ref{s:prelim}.

In Section  \ref{ss:recognition}, we show how a subprojection of $\cS$ with the form $\TMAX(\cE)$ can be recognized inside~$\cS$. A~further reformulation is carried out in Section \ref{ss:inductive}, where we show that if sequences $(E_n)_{n\geq 1}$ in $\Om$ and $(\cG_n)_{n\geq1}$ in $\cP(\Om)$ can be constructed to satisfy certain recursive properties, then one can inductively build a subprojection of $\cS$ of the form~$\TMAX(\cE)$.
Similar results are established for $\TMIN$ and~$\TORT$;
it is notable that for these two cases, repeated subprojections of $\cS$ occur naturally. (The relation $\subproj$ is transitive, as will be shown in Lemma \ref{l:subproj is transitive}.)

We can view the task of constructing suitable sequences (yielding either $\TMAX(\cE)$, $\TMIN(\cE)$ or $\TORT(\cE)$) as a $1$-player game in discrete time, where the ``state space'' of the game consists of subprojections of $\cS$.
We start at time $n=1$ with initial state $\cS_1=\cS$; and at stage $n$ of the game, if we are in state $\cS_n$, we pass to a subprojection $\cS_{n+1}\subproj\cS_n$, subject to certain rules.
At stage $n$, the permissible ``moves'' have three possible types, labelled $\MAX_n$, $\MIN_n$ and $\ORT_n$.

The results of Section \ref{ss:inductive} have the following interpretation: we win the game --- in the sense of being able to construct $\cT_*$ of the desired form --- if it is possible from some point onwards to make an infinite sequence of consecutive moves all of the same type, e.g.
\[ \MAX_1 \,,\, \MIN_2 \,,\, \ORT_3 \,,\, \ORT_4 \,,\, \MAX_5 \,,\, \MAX_6 \,,\, \MAX_7 \,, \dots \]
In making such moves, we might run into a dead end and need to backtrack to an earlier stage of the game. To make progress we need a closer look at how this problem could arise, and what conditions can be imposed to rule out such problems.

Let $\cS'\subproj\cS$.
In Definition \ref{d:automatic}, we introduce the concept of being \dt{$\TMAX$-automatic}. This definition is rather technical, but is guided by the results in Section \ref{ss:inductive}; the key point is that if $\cS'$ is $\TMAX$-automatic, then for any $n\in\Nat$ and any $\cV\subproj\cS'$, there is \emph{always} some move of type $\MAX_n$ starting from $\cV$ (Proposition \ref{p:auto-generate}). Since such a move takes us from $\cV$ to some $\cV'\subproj \cV\subproj \cS'$, and since $\subproj$ is transitive, one obtains a winning sequence of moves without needing to look ahead; at stage $n$ we can always make a move of $\MAX_n$, and whichever one we pick, there will be a move of type $\MAX_{n+1}$ available at the next stage.

Similarly: there are notions of being \dt{$\TMIN$-automatic} and \dt{$\TORT$-automatic} (Definition~\ref{d:automatic} again), and analogous results showing that set systems with these properties allow the construction of copies of $\TMIN$ or $\TORT$ by ``making moves without needing to look ahead'' (Proposition \ref{p:auto-generate} again).
Therefore, Theorem \ref{t:headline result} will follow if we can show that  $\cS$ always has a subprojection that is either $\TMAX$-automatic, $\TMIN$-automatic or $\TORT$-automatic. This will be achieved in Proposition \ref{p:auto-inevitable}; it is here that the precise technical definitions of $\cT_*$-automaticity (for $*\in \{{\rm MAX}, {\rm MIN}, {\rm ORT}\}$) becomes important. Section~\ref{section:lemmas} is devoted to the technical lemmas needed to make this work.

Note that we will end up with a sharper version of Theorem \ref{t:headline result}. It is easy to construct examples of $\cS$ which have subprojections of the form $\TMAX(\cE)$ but are not $\TMAX$-automatic; nevertheless, Proposition \ref{p:auto-inevitable} tells us that \emph{some subprojection} of $\cS$ will be $\TMAX$-automatic.

\end{subsection}

\end{section}

\begin{section}{Preliminary results}\label{s:prelim}

\begin{subsection}{Some notation and terminology}
In this subsection, we fix the notational conventions that will be used in the rest of the paper.

\begin{notn}
Let $\cS$ be a set system on a nonempty set $\Omega$. We use the following
terminology and notational conventions to aid clarity of our presentation.
\begin{romnum}
\item Elements of $\Omega$ will usually be denoted by lower-case Greek letters. General subsets of $\Omega$ may be denoted by upper-case Roman letters or upper-case Greek letters; but see the next item.
\item
By \dt{members} of $\cS$, we mean subsets of $\Omega$ which belong to $\cS$. We use letters such as $\ssa$, $\ssb$, $\ssp$, etc.,  to denote members of $\cS$, and write $\ssa\cup\ssb$ and $\ssa\cap \ssb$ for their union and intersection respectively. The complement of $\ssa$ in $\Om$ is denoted by $\ssa^c$.
If it happens that $\ssa$ and $\ssb$ are disjoint subsets of $\Omega$, we shall sometimes emphasize this by writing their union as $\ssa\dotcup\ssb$. 

\item
We recall (Definition~\ref{d:union-closed}) that $\cS$ is said to be
  \dt{union-closed} if $\ssa,\ssb\in\cS \implies \ssa\setcup\ssb\in\cS$.
Note that with this terminology, ``union-closed'' need not imply being closed under taking arbitrary unions: the set system $\FIN(\Omega)$ in Example \ref{eg:fin-and-cofin} illustrates this point.
\item
Let $\cF\subseteq\cP(\Omega)$. We write $\join(\cF)$ for $\bigsetcup_{\ssx\in\cF}\ssx$ and $\meet(\cF)$ for $\bigsetcap_{\ssx\in\cF}\ssx$. Note that if
$\cS$ is union-closed and $\cF\subseteq\cS$ is \emph{finite},
then $\join(\cF)\in\cS$.

\item
Let $(\cS_i)$ be a family of subsystems of $\cP(\Om)$. We denote the union of these subsystems by $\bigsscup_i \cS_i$ ; this should not be confused with $\{ \bigsetcup_i \ssx_i \colon \ssx_i\in\cS_i \;\forall\,i\}$. The intersection of the $\cS_i$ is denoted similarly by $\bigsscap_i \cS_i$. 
\end{romnum}

\end{notn}

Given $\cS\subseteq\cP(\Om)$ and $\Gm\subseteq\Om$, we already defined the \dt{projection} $\cS\owedge\Gm$ in Definition \ref{d:subprojection}. We may also define other set systems derived from $\cS$ using $\Gm$:

\begin{equation}\label{eq:derived set systems}
\begin{aligned}
\cS\ominus \Gm &\defeq \{ \ssx\setminus \Gm \st \ssx\in \cS\} &\subseteq & \cP(\Om\setminus\Gm),
\\
\cS\ovee \Gm &\defeq \{ \ssx\setcup \Gm \st \ssx\in \cS\} &\subseteq & \cP(\Om),
\\
\cS^{-\Gm} & \defeq  \{\ssx\in\cS\colon\  \ssx\cap\Gm=\emptyset\} & \subseteq & \cS,  \\
\cS_\Gm &\defeq \{\ssx\in\cS\colon\  \ssx\supseteq\Gm\} &\subseteq & \cS.
\end{aligned}
\end{equation}
Clearly, each of these set systems is union-closed if $\cS$ is.

\end{subsection}

\begin{subsection}{More on subprojections}
\label{ss:more subproj}

\begin{lem}[Equivalent notions of subprojection]
\label{l:subprojection}
Let $\cS,\cT\subseteq\cP(\Om)$ and let $X\subseteq\Om$. The following are equivalent:
\begin{romnum}
\item\label{li:sub-of-proj}
(subsystem of a projection) there exists $X\subseteq \Omega$ such that $\cT\subseteq \cS\owedge X$;
\item\label{li:proj-of-sub}
(projection of a subsystem) there exist  $X\subseteq \Omega$ and $\cS'\subseteq\cS$ such that $\cT=\cS'\owedge X$, i.e. $\cT\subproj\cS$.
\item\label{l:translation}
$\cT \subseteq \cS\owedge\join(\cT)$.
\end{romnum}
Moreover, if both $\cS$ and $\cT$ are union-closed, in \ref{li:proj-of-sub} we may choose $\cS'$ to be union-closed.
\end{lem}

The proof of the lemma is routine; we provide details in order to familiarize the reader with the notation that has been previously introduced.

\begin{proof}
Suppose that \ref{li:proj-of-sub} holds. Then
\[ \cT=\cS'\owedge X = \{ \ssa\setcap X \colon \ssa\in \cS'\} \subseteq \{ \ssa\setcap X \colon \ssa\in \cS \} = \cS\owedge X \;,\]
and thus \ref{li:sub-of-proj} holds.
Similarly: if \ref{li:sub-of-proj} holds for some $X\subseteq\Om$, then since $\join(\cT)\subseteq \join(\cS\owedge X)\subseteq X$ we have
 \[
\cT = \cT\owedge\join(\cT) \subseteq ( \cS\owedge X)\owedge \join(\cT)= \cS\owedge\join(\cT)\;; \]
and thus \ref{l:translation} holds. 
The converse implication \ref{l:translation}$\implies$\ref{li:sub-of-proj}  is trivial. 

Finally, suppose that \ref{li:sub-of-proj} holds. Define
$\cS' \defeq \{ \ssa\in \cS \colon \ssa\setcap X \in \cT\} \subseteq \cS$,
noting that if $\cT$ and $\cS$ are union-closed then so is $\cS'$.
By construction,
\[
\cS'\owedge X
 = \{ \ssa\setcap X \colon \ssa\in\cS\text{ and } \ssa\setcap X \in\cT\} \subseteq \cT \;.
\]
But by \ref{li:sub-of-proj}, each $\ssb\in\cT$ is equal to $\ssa\setcap X$ for some $\ssa\in\cS$; this $\ssa$ then belongs to $\cS'$, and so $\ssb\in \cS'\owedge X$. Thus $\cT\subseteq\cS'\owedge X$, and so \ref{li:proj-of-sub} holds.
\end{proof}

The property in Lemma \ref{l:subprojection}\ref{li:sub-of-proj} is often more convenient to work with than the original definition of a subprojection,  as is shown clearly in the short proof of the next important lemma.

\begin{lem}[Transitivity of $\subproj$]
\label{l:subproj is transitive}
Let $\cT_1,\cT_2,\cT_3$ be set systems on $\Om$. If $\cT_1\subproj\cT_2$ and  $\cT_2\subproj\cT_3$, then $\cT_1\subproj\cT_3$.
\end{lem}

\begin{proof}
By Lemma \ref{l:subprojection} \ref{li:sub-of-proj}, there exist $X_2\subseteq \Omega$ such that $\cT_1\subseteq \cT_2\owedge X_2$ and $X_3\subseteq \Omega$ such that $\cT_2\subseteq \cT_3\owedge X_3$. Thus
\[
\cT_1\subseteq (\cT_3\owedge X_3)\owedge X_2=\cT_3\owedge(X_3\cap X_2).
\]
Hence, using the same lemma, $\cT_1\subproj \cT_3$.
\end{proof}

We chose Lemma \ref{l:subprojection}\ref{li:proj-of-sub} as our definition of subprojection because, later in the paper, several arguments have an inductive step with the following form: starting from some set system $\cT\subseteq\cP(\Om)$, we form $\cT^{-F}\owedge X$ for certain choices of $F$ and $X$. In other words, the natural order of operations is to pass to a subsystem of $\cT$ and then apply a projection.

By Lemma \ref{l:subprojection},
$\cT^{-F}\owedge X$ is a subsystem of a projection of~$\cT$,
but it is important to note that in general it \emph{is not the same} as $(\cT\owedge X)^{-F}$.
However, in certain cases the two set systems coincide. The next lemma provides sufficient conditions for this.

\begin{lem}[Swapping projection and exclusion/inclusion]
\label{l:swap-proj}
Let $X,F\subseteq\Om$ and let $\cT\subseteq\cP(\Om)$.
\begin{romnum}
\item\label{li:swap exclude}
$(\cT\owedge X)^{-F} \supseteq \cT^{-F} \owedge X$. If $X\supseteq F$, then equality holds.
\item\label{li:swap include}
$(\cT\owedge X)_{F} \subseteq \cT_{F} \owedge X$. If $X\supseteq F$, then equality holds.
\end{romnum}
\end{lem}

\begin{proof}
We have
\[ \begin{aligned}
(\cT\owedge X)^{-F}
& =\{ \ssa\setcap X \colon \ssa\in \cT, \ssa\setcap X \setcap F =\emptyset \}  \\
& \supseteq  \{ \ssa\setcap X \colon \ssa\in \cT, \ssa\setcap F =\emptyset \} 
& =\cT^{-F} \owedge X \;.
\end{aligned} \]
If $X\supseteq F$ then $X\setcap F = F$ and so the inclusion above is an equality. This proves part \ref{li:swap exclude}.

Similarly,
\[ \begin{aligned}
(\cT\owedge X)_{F}
& =\{ \ssa\setcap X \colon \ssa\in \cT, \ssa\setcap X \supseteq F \}  \\
& \subseteq  \{ \ssa\setcap X \colon \ssa\in \cT, \ssa\supseteq F \} 
& =\cT_{F} \owedge X \;.
\end{aligned} \]
If $X\supseteq F$ then any $\ssa\in\cT$ satisfying $\ssa\supseteq F$ also satisfies $\ssa\setcap X\supseteq F$, and so the inclusion above is an equality. This proves part \ref{li:swap include}.
\end{proof}

\end{subsection}

\begin{subsection}{Incompressible subsets and their witnesses}
\label{ss:witnesses}

\begin{dfn}[Incompressible subsets of $\cP(\Om)$]
\label{d:incompressible}
Let $\cF$ be a finite subset of $\cP(\Om)$. If there is a proper subset $\cF'\subset \cF$ such that $\join(\cF')=\join(\cF)$, we say that $\cF$ is \dt{compressible}; otherwise, we say $\cF$ is \dt{incompressible}.

Given $\cS\subseteq\cP(\Omega)$,
we write $\Inc(\cS)$ for the set of all finite incompressible subsets of $\cS$. Given $n\in{\mathbb N}$, we write $\Inc_n(\cS)=\{ \cF \in \Inc(\cS) \colon |\cF|=n\}$.
\end{dfn}

Thus, if $\cS$ is a union-closed set system on $\Om$, Definition~\ref{d:breadth} may be rewritten as
$\br(\cS)= \sup\{ n\in \Nat \colon \Inc_n(\cS)\text{ is non-empty}\}$.

\begin{rem}
\label{r:nested}
The example $\FIN(\Om)$ (Example~\ref{eg:fin-and-cofin}) has the property that it contains an increasing sequence  $\cF_1\subset\cF_2\subset\dots$ of finite incompressible subsets. This is not typical; such sequences do not exist for $\TMAX(\cE)$, $\TMIN(\cE)$ or $\TORT(\cE)$.
\end{rem}

If $\cS$ is union-closed, incompressible subsets of $\cS$ with size $m$ generate subsets of $\cS$ that resemble $\cP(\{1,\dots,m\})\setminus\{\emptyset\}$. We now make this more precise.
Let $\cF\subseteq\cP(\Om)\setminus\{\emptyset\}$
with $|\cF|=m\geq 2$, and enumerate $\cF$ as $\{\ssx_1,\dots, \ssx_m\}$.
It~is easily checked that the following statements are equivalent:
\begin{romnum}
\item there exists $j$ such that $\ssx_j\setcap \bigsetcap_{i\neq j} {\ssx_i}^c=\emptyset$\/;
\item there exists $j$ such that $\ssx_j\subseteq \bigsetcup_{i\neq j} \ssx_i$\/;
\item $\{\ssx_1,\dots, \ssx_m\}$ is compressible.
\end{romnum}
Thus, such an $\cF$ is incompressible if and only if there exist $\gm_1,\dots, \gm_m\in\Om$ such that
\begin{equation}\label{eq:witness}
\gamma_i\in \ssx_i\setminus \ssx_j \qquad\text{whenever $i\neq j$.}
\end{equation}
In fact, this characterizes incompressible subsets of $\cP(\Om)$ which have size $\geq 2$ (for if $\cG\in\Inc(\Om)$ with $|\cG|\geq 2$, then $\emptyset\notin\cG$).

\begin{dfn}\label{d:witness}
Let $\cF=\{\ssx_1,\ldots, \ssx_m\}\subset \cP(\Om)$ with $m\geq 2$ be incompressible. A subset $F=\{\gamma_1,\ldots,\gamma_m\}$ of $\Om$ is called a \dt{witness of incompressibility} for $\cF$, or just a \dt{witness for~$\cF$}, if it satisfies (\ref{eq:witness}).
Equivalently, each $\gm_j$ is contained in a \emph{unique} member of $\cF$.
\end{dfn}

\begin{rem}
Every $1$-element subsystem of $\cP(\Om)$ is incompressible. However,  for the definition of a witness we restrict our attention to subfamilies of size at least~$2$.
One reason is that some useful results such as Lemma \ref{l:easy witness}\ref{li:witness-join-meet} become false if we allow witnesses of size~$1$.  Also, the 1-element subfamily $\{\emptyset\}$ is incompressible, but does not admit a witness at all. 

\end{rem}

The next remark is included for background interest; we will not need it for the proof of our main result. It provides additional detail connected with the remarks in Section~\ref{ss:wider context}.

\begin{rem}[Breadth and VC-dimension]\label{r:breadth and VC}
If $\cS$ is a set system on $\Om$ (not necessarily union-closed) and $F$ is a finite non-empty subset of $\Om$, then $\cS$ is said to \dt{shatter~$F$} if $\cS\owedge F = \cP(F)$. This concept is important for statistical learning theory, so it is interesting to note its connection with the (lattice-theoretic) notion of breadth.

Let $\cS\subseteq\Om$ be union-closed. If $\cS$ shatters some finite set $\emptyset\neq F\subseteq\Om$, then {\itshape a~fortiori} $\cS\owedge F \supseteq \{\{\gm\}\colon \gm\in F\}$, and hence $F$ is a witness of incompressibility for some $\cF\subseteq\cS$.
In the other direction: let $\cS\subseteq\cP(\Om)$ be union-closed and satisfy $\emptyset\in\cS$. Suppose there exists $\cG\in\Inc_n(\cS)$ for some $n\geq 2$, and let $G$ be a witness of incompressibility for $\cG$. It then follows from the assumptions on $\cS$ that
$\cS$ shatters $G$.

The \dt{Vapnik--Chervonenkis (VC) dimension} of a set system $\cS\subseteq\cP(\Om)$ is defined as the supremum of $|F|$ over all finite subsets $F\subseteq\Om$ that are shattered by $\cS$. The calculations above show that when $\cS$ is union-closed and $\emptyset\in\cS$, then $\br(\cS)=\dim_{\rm VC}(\cS)$.
\end{rem}
\end{subsection}

\begin{subsection}{Various lemmas}
In this section,
we collect some properties of witnesses which will be useful later,
and then apply them to study how breadth behaves when combining or projecting set systems.
\begin{lem}\label{l:easy witness}
Let $\cF\in\Inc_n(\cP(\Om))$ for some $n\geq 2$.
\begin{romnum}
\item\label{li:witness-join-meet}
Any witness for $\cF$ is contained in $\join(\cF)$ and has
empty intersection with $\meet(\cF)$.
\item\label{li:subset of witness}
If $F$ is a witness for $\cF$, and $G \subseteq F$ with $|G|\ge 2$, then $G$ is a witness for a (unique) incompressible $\cG\subseteq\cF$.
\end{romnum}
\end{lem}

\begin{proof}
Let $F$ be a witness for $\cF$. Each element of $F$ is an element of some member of $\cF$ and hence is an element of the union of all members of $\cF$; that is, $F\subseteq\join(\cF)$. If there existed some $\gm\in F \cap \meet(\cF)$ then $\gm$ would be an element of two distinct members of $\cF$, contradicting what it means to be a witness. This proves~\ref{li:witness-join-meet}.

Now suppose that $G\subseteq F$ with $|G|\geq 2$. For each $\alpha\in G$ let $\ssx_\alpha$ be the \emph{unique} member of $\cF$ to which $\alpha$ belongs, and let $\cG=\{\ssx_\alpha \colon \alpha \in G\}\subseteq\cF$. It is easily checked that $G$ is a witness for $\cG$. If $\cH\subseteq\cF$ also has the property that $G$ is a witness for $\cH$, then $\cH$ must contain $\cG$ (since $\ssz\in \cF$ and $\alpha\in \ssz$ implies $\ssz=\ssx_\alpha$) and must also satisfy $|\cH|= |G| = |\cG|<\infty$, and so $\cH=\cG$. This proves \ref{li:subset of witness}.
\end{proof}

Whenever $\cT\subproj \cS$, incompressible subsets of $\cT$ can always be ``lifted'' to incompressible subsets of~$\cS$. To state things precisely, we make the following definition.

\begin{dfn}[Liftings]\label{d:liftings}
Let $\cS\subseteq\cP(\Om)$ be a non-empty set system, and let $\Gamma\subseteq \Om$. Given $\cJ\subseteq \cS\owedge \Gamma$, a \dt{lifting of $\cJ$} (relative to $\cS$ and $\Gamma$) is a set system $\cI \subseteq \cS$ such that the truncation map $\ssa\mapsto \ssa\setcap \Gamma$ restricts to a bijection $\cI\to\cJ$.
\end{dfn}

Liftings of a given $\cJ\subseteq\cS\owedge \Gamma$ always exist, since the truncation map $\cS\to\cS\owedge \Gamma$ is surjective by definition. (In the desired applications, $\cJ$ will be finite, so we can avoid invoking the axiom of choice.)

\begin{lem}[Lifting incompressible subsets of a subprojection]\label{l:lift inc}
Let $\cS\subseteq\cP(\Om)$, let $\Gamma\subseteq\Om$ be non-empty, and let $\cT=\cS\owedge \Gamma$.
Suppose $\cJ\in\Inc(\cT)$ with $|\cJ|\geq 2$, and let $J$ be a witness for $\cJ$. Let $\cI\subseteq \cS$ be any lifting of $\cJ$. Then $\cI$ is an incompressible subsystem of $\cS$, and  $J$ is a witness for $\cI$.
\end{lem}

\begin{proof}
By Lemma \ref{l:easy witness}\ref{li:witness-join-meet}, we have $J\subseteq \join(\cJ)\subseteq  \Gamma$.
In particular, $J\subseteq \join(\cI)$. It suffices to prove that each $\alpha\in J$ belongs to a unique member of $\cI$.
Suppose otherwise: then there exists $\alpha\in J$ which belongs to two distinct members of $\cI$, say $\ssx_1$ and $\ssx_2$.
Hence $\alpha\in \ssx_1\cap\Gamma$ and $\alpha\in \ssx_2\cap\Gamma$. This is a contradiction, since $\ssx_1\cap \Gamma$ and $\ssx_2\cap \Gamma$ are members of $\cJ$ and $J$ is a witness for $\cJ$.
\end{proof}

\begin{lem}[Coarse-graining of an incompressible system]
\label{l:coarse-grain}
Let $\cG$ be an incompressible subsystem of $\cP(\Om)$, with a witness $G$. 
Let $\cG= \bigsscup_{i=1}^n \cG_i$, with $n\geq 2$, be a partition of $\cG$ into a disjoint union of non-empty subsystems.
For each $i=1,\ldots, n$,  put $\ssx_i = \join(\cG_i)$, and pick exactly one element $\gamma_i \in \ssx_i\cap G$. Then $\{ \ssx_i \colon i=1,\dots, n\}$ is incompressible, with $\{\gm_i \colon 1\leq i\leq n\}$ as a witness.
\end{lem}

\begin{proof}
By definition, for each $i\in\{1,\dots, n\}$, $\gm_i\in \ssx_i$.
Suppose that there exist distinct $i$ and $j$ with $\gm_i \in \ssx_i\cap\ssx_j$. Then there exist $\ssa_i\in\cG_i\subset \cG$ and $\ssa_j\in \cG_j\subset \cG$ such that $\gm_i \in \ssa_i \setcap \ssa_j$, which contradicts the assumption that $G$ is a witness for $\cG$. We conclude that whenever $i$ and $j$ are distinct, $\gm_i$ belongs to $\ssx_i\setminus\ssx_j$; thus it is a witness of incompressibility for $\{ \ssx_1,\dots,\ssx_n\}$. 
\end{proof}

In the last lemma of this section, we observe how the breadth of a union-closed set system interacts with the breadth of its subsystems and projections. 

\begin{lem} \label{l:swissarmyknife}
Let $\cS, \cS_1,\dots, \cS_m  \subseteq \cP(\Om)$ be union-closed,
and let $X_1,\dots, X_m$ be subsets of $\Om$ such that $\Om = X_1\setcup \dots \setcup X_m$.
\begin{romnum}
\item\label{li:breadth of merged}
Suppose $\br(\cS_i)<\infty$ for each $i$. If $\cS=\bigsscup_{i=1}^m \cS_i$, then $\br(\cS)<\infty$.
\item\label{li:breadth of partition}
Suppose that $\br(\cS\owedge X_i)<\infty$ for all $1\leq i\leq m$. Then $\br(\cS)<\infty$.
\end{romnum}
\end{lem}

\begin{proof}
%
To prove \ref{li:breadth of merged}, let $n_i=\br(\cS_i)$ and let $n=\sum_{i=1}^m n_i$. For $\cJ\subseteq \cS$, we may write $\cJ=\bigsscup_{i=1}^m (\cJ \sscap \cS_i)$. By the pigeonhole principle, if  $|\cJ| >n$ then there is at least one $i$ such that $|\cJ \sscap \cS_i | > n_i$, and for this $i$, $\cJ\sscap\cS_i$ is compressible. Hence $\cJ$ is compressible, as any superset of a compressible set is again compressible. This implies that $\br(\cS)\leq n$, so $\cS$ has finite breadth. 

We now prove \ref{li:breadth of partition}. Towards a contradiction, suppose that $\br(\cS)=\infty$. Let $n_i=\br(\cS\owedge X_i)$, and take $k>\sum_{i=1}^m n_i$.
Since $\cS$ has infinite breadth, it contains an incompressible subset $\cF$ with a witness $F$ such that  $|\cF|=|F|=k$.

By the definition of a witness, for each $\alpha\in F$ there exists a unique member of $\cF$ containing $\alpha$, which we denote by~$\ssx_\alpha$.
For $1\leq i\leq m$, define
$\cF_i=\{\ssx_\alpha \colon \alpha \in F\cap X_i\}\subseteq \cF$.
Then $\cF = \{ \ssx_\alpha \colon \alpha\in F\} = \bigsscup_{i=1}^m \cF_i$,
as $\Om = X_1\setcup \dots \setcup X_m$. So by the pigeonhole principle, there is at least one $i$ such that $|\cF_i|>n_i$. 
Hence $|F\cap X_i|=|\cF_i| \geq n_i+1\geq 2$. Then since $F\cap X_i$ is a witness of incompressibility for $\cF_i\owedge X_i$, we have $n_i=\br(\cS\owedge X_i) \geq |\cF_i \owedge X_i| \geq n_i+1$;
a contradiction.
\end{proof}

\end{subsection}
\end{section}

\begin{section}{A closer look at our target configurations}

\begin{subsection}{Recognizing the desired configurations}
\label{ss:recognition}
Throughout this subsection, \emph{but not the subsequent ones}, we fix:
\begin{itemize}
\item a union-closed set system $\cS\subseteq\cP(\Om)$;
\item a sequence $(E_n)_{n\geq 1}$ of finite, non-empty subsets of $\Om$
with $|E_n|=n+1$ for all $n\in\Nat$.
\end{itemize}
For $n\in\Nat$, we put $E_{<n}\defeq\bigsetcup_{j=1}^{n-1} E_j$ and $E_{>n}\defeq\bigsetcup_{k=n+1}^\infty E_k$, with the convention that $E_{<1}=\emptyset$.
Let $\cE=\{ E_n \colon n\in\Nat\}$, so that $\join(\cE)=\bigsetcup_{n\geq 1} E_n$.

Unlike Definition~\ref{d:3kings}, we do not assume at the outset that the sets $E_n$ are pairwise disjoint.
However, the results of this subsection, Lemmas \ref{l:TMAX-orig}, \ref{l:TMIN-orig} and \ref{l:TORT-orig}, provide sufficient\footnotemark\ conditions for the $E_n$ to be pairwise disjoint and for $\cS$ to have a subprojection equal to $\TMAX(\cE)$, $\TMIN(\cE)$ and $\TORT(\cE)$ respectively.
\footnotetext{It will be shown {\itshape en route} that these conditions are also necessary, although this is not required for the proof of Theorem~\ref{t:headline result}.}


\paragraph{Recognizing $\TMAX(\cE)\subproj\cS$.}
Suppose that the sets $E_n$ are pairwise disjoint, and that $\TMAX(\cE)\subseteq\cS\owedge(\join(\cE))$.
The $n$th level in $\TMAX(\cE)$ is generated by taking finite unions of all sets of the form $\{\gamma\}\dotcup E_{<n}$ where $\gm\in E_n$; and such a set must have the form $\ssx\setcap \join(\cE)$ for some $\ssx\in\cS$. 
Collecting such $\ssx$, one for each $\{\gamma\}\dotcup E_{<n}$,  gives rise to a lifting for the generating set with witness $E_n$.
Hence, there exists a sequence $(\cF_n)_{n\geq 1}$ of subsets of $\cS$ with the following properties:
\begin{quote}
\begin{itemize}
\item[(WIT)] each $\cF_n$ is incompressible, and has $E_n$ as a witness;
\item[(MAX)] for all $n\in \Nat$, $\meet(\cF_n)\supseteq E_{<n}$ and $\join(\cF_n)\cap E_{>n}=\emptyset$.
\end{itemize}
\end{quote}
Conversely, the next lemma says in effect that these two necessary conditions are also sufficient. 
In particular, they automatically imply that $E_j\setcap E_k=\emptyset$ for all $1\leq j < k$.

\begin{lem}\label{l:TMAX-orig}
Suppose there exists a sequence $(\cF_n)_{n\geq 1}$ of subsets of $\cS$ satisfying (WIT) and (MAX).
Then the sets $E_n$ are pairwise disjoint; and the union-closed set system generated by $\left(\bigsscup_{n=1}^\infty\cF_n\right)\owedge \join(\cE)$ is $\TMAX(\cE)$.
\end{lem}

\begin{proof}
From (MAX) we have $E_{< n} \subseteq \meet(\cF_n)$ for every $n\geq 2$. But by (WIT),  $E_n$ is a witness for $\cF_n$, so Lemma \ref{l:easy witness} implies that $E_n\setcap E_{<n} \subseteq E_n\setcap \meet(\cF_n)=\emptyset$. Hence $E_j\setcap E_k=\emptyset$ whenever $1\leq j < k$.

To prove the second part, given any $\ssx\in\cF_n$ consider the  decomposition
\[
\ssx \setcap\join(\cE) = (\ssx\setcap E_{<n} )\dotcup (\ssx\setcap E_n) \dotcup (\ssx\setcap E_{>n}).
\]
(Recall that we use $\dotcup$ to signify the fact that the sets forming the union are pairwise disjoint.)
The second term equals $\{\gamma\}$ for some $\gamma\in E_n$, since $\ssx\in\cF_n$ and $E_n$ is a witness for $\cF_n$. By the condition (MAX), the first term equals $E_{<n}$ and the third term is empty. Thus each member of $\cF_n\owedge\join(\cE)$ is one of the generating elements of $\TMAX(\cE)$.
\end{proof}

\paragraph{Recognizing $\TMIN(\cE)\subproj \cS$.}
Suppose that the sets $E_n$ are pairwise disjoint, and that $\TMIN(\cE)\subseteq\cS\owedge(\join(\cE))$.
 The $n$th level in $\TMIN(\cE)$ is generated by taking finite unions of
all sets of the form $\{\gm\}\dotcup E_{>n}$ where $\gm\in E_n$, and such a set must have the form $\ssx\cap \join(\cE)$ for some $\ssx\in \cS$. Collecting such $\ssx$, one for each $\{\gamma\}\dotcup E_{>n}$,  gives rise to a lifting  for the generating set with witness $E_n$.
Hence, there exists a sequence $(\cF_n)_{n\geq 1}$ of subsets of $\cS$ with the following properties:
\begin{quote}
\begin{itemize}
\item[(WIT)] each $\cF_n$ is incompressible, and has $E_n$ as a witness;
\item[(MIN)] for all $n\in \Nat$, $\meet(\cF_n)\supseteq E_{>n}$ and $\join(\cF_n)\cap E_{<n}=\emptyset$.
\end{itemize}
\end{quote}
There is an analogue of Lemma \ref{l:TMAX-orig}. The proof is very similar, so we will give fewer details.
\begin{lem}\label{l:TMIN-orig}
Suppose there exists a sequence $(\cF_n)_{n\geq 1}$ of subsets of $\cS$ satisfying (WIT) and (MIN).
Then the sets $E_n$ are pairwise disjoint; and the union-closed set system generated by $\left(\bigsscup_{n=1}^\infty\cF_n\right)\owedge \join(\cE)$ is $\TMIN(\cE)$.
\end{lem}

\begin{proof}
By (WIT) $E_n$ is a witness for $\cF_n$, so in particular $E_n\subseteq\join(\cF_n)$. But by the second part of (MIN), $\join(\cF_n)$ is disjoint from $E_{<n}$; thus $E_{<n}\setcap E_n=\emptyset$, as in the $\TMAX$ case.

For the second part, note that for each $\ssx\in\cF_n$ we have
\[
\ssx\setcap \join(\cE) = (\ssx\setcap E_{<n}) \dotcup (\ssx\setcap E_n) \dotcup (\ssx \setcap E_{>n}) = \{\gm\}\dotcup E_{>n}
\]
and these sets generate $\TMIN(\cE)$.
\end{proof}

\paragraph{Recognizing $\TORT(\cE)\subproj\cS$.}
Suppose that the sets $E_n$ are pairwise disjoint, and that $\TORT(\cE)\subseteq\cS\owedge(\join(\cE))$.
The $n$th level in $\TORT(\cE)$ is generated by taking finite unions of all sets of the form $\{\gm\}\dotcup E_{<n}\dotcup E_{>n}$ where $\gm\in E_n$, and such a set must have the form $\ssx\cap \join(\cE)$ for some $\ssx\in \cS$.
Collecting such $\ssx$, one for each $\{\gm\}\dotcup E_{<n}\dotcup E_{>n}$,  gives rise to a lifting for the generating set with witness $E_n$.
Hence, there exists a sequence $(\cF_n)_{n\geq 1}$ of subsets of~$\cS$ with the following properties:
\begin{quote}
\begin{itemize}
\item[(WIT)] each $\cF_n$ is incompressible, and has $E_n$ as a witness;
\item[(ORT)] for all $n\in \Nat$, $\meet(\cF_n)\supseteq E_{>n}\dotcup E_{<n}$.
\end{itemize}
\end{quote}
As in the previous two cases, we get a converse statement.
\begin{lem}\label{l:TORT-orig}
uppose there exists a sequence $(\cF_n)_{n\geq 1}$ of subsets of $\cS$ satisfying (WIT) and (ORT). 
Then the sets $E_n$ are pairwise disjoint; and the union-closed set system generated by $\left(\bigsscup_{n=1}^\infty\cF_n\right)\owedge \join(\cE)$ is $\TORT(\cE)$.
\end{lem}

\begin{proof}
Similarly to the $\TMAX$ case, by combining the condition $E_{>n} \dotcup E_{<n} \subseteq \meet(\cF_n)$ with the condition (WIT), it follows that $E_n\setcap (E_{>n}\cup E_{<n})=\emptyset$. We also have, for each $n\in\Nat$ and each $\ssx\in\cF_n$,
\[ \ssx \setcap \join(\cE) = (\ssx\setcap E_{<n}) \dotcup (\ssx\setcap E_n) \dotcup (\ssx\setcap E_{>n}) =  E_{<n} \dotcup \{\gamma\} \dotcup E_{>n}\]
and sets of these form generate $\TORT(\cE)$.
\end{proof}

\end{subsection}

\begin{subsection}{Inductive reformulation of our target configurations}\label{ss:inductive}
We now build on the previous observations to obtain necessary and sufficient conditions with an inductive flavour. Intuitively, these conditions have the following flavour: in each of the three cases there is a certain kind of ``move'' that one is allowed to apply to a given union-closed set system; and being able to make an infinite sequence of ``moves'' of that kind will allow one to construct a subprojection of the form $\TMAX$, $\TMIN$ or $\TORT$.

Let $\Om$ be a non-empty set. Recall (Example~\ref{eg:fin-and-cofin}) that $\FIN(\Om)$ denotes the set of all finite subsets of $\Om$.

\begin{prop}[Reformulating $\TMAX$]\label{p:reformulate TMAX}
Let $\cS\subseteq \cP(\Om)$ be a union-closed set system, and let $(E_n)_{n\geq 1}$ be a sequence in $\FIN(\Om)$, with  $|E_n|\geq 2$ for all $n$. The following are equivalent:
\begin{romnum}
\item\label{li:TMAX-ish}
there exist incompressible subsets $\cF_n\subseteq \cS$ for all $n\in \Nat$, which satisfy the conditions (WIT) and (MAX) with respect to $(E_n)_{n\ge 1}$;
\item\label{li:TMAX-ind-new}
there exist incompressible subsets $\cG_n\subseteq \cP(\Om)$, for $n\geq 1$, with $E_n$ a witness for $\cG_n$, such that $\cG_1\subseteq \cS$ and 
\[
\cG_{n+1}\subseteq \cS\ominus\left(\bigcup_{j=1}^n\join(\cG_j)\right) \qquad \text{ for all } n\in \Nat.
\]
\end{romnum}
\end{prop}
\begin{proof}
\ref{li:TMAX-ish}$\implies$\ref{li:TMAX-ind-new}: 
Put $\cG_1 \defeq \cF_1\subseteq\cS$ and
$\cG_n \defeq \cF_n\ominus \bigsetcup_{j=1}^{n-1} \join(\cF_j)$
for all $n\geq 2$.
Since $E_n$ is a witness for $\cF_n$ and is disjoint from $\join(\cF_{j})$ for every $j<n$ by condition~(MAX), we see that $E_n$~remains a witness
 for~$\cG_n$. Moreover, by induction we have $\bigsetcup_{j=1}^n \join(\cG_j)=\bigsetcup_{j=1}^n\join(\cF_j)$ for all $n\geq 1$. Hence
\[
\cG_{n+1}
= \cF_{n+1}\ominus \bigsetcup_{j=1}^n \join(\cG_j)
\subseteq \cS\ominus \bigsetcup_{j=1}^n \join(\cG_j)
\]
as required.
	
\medskip
\ref{li:TMAX-ind-new}$\implies$\ref{li:TMAX-ish}
To ease notation, define an increasing sequence $\ssb_n = \bigcup_{j=1}^n \join(\cG_j)$ in $\cP(\Om)$, with the convention that $\ssb_0=\emptyset$. Then our assumptions can be rephrased as:  $\cG_n\subseteq \cS\ominus \ssb_{n-1}$ for all $n\geq 1$. 
Also, note that $E_n \subseteq \join(\cG_n) = \ssb_n \setminus \ssb_{n-1}$ for all $n\in\Nat$.

We shall prove by induction that
$\cG_n \ovee \ssb_{n-1} \subseteq \cS$ for all $n\in\Nat$. For $n=1$ this holds since $\ssb_0=\emptyset$ and $\cG_1\subseteq\cS$. 
Suppose for some $n\in\Nat$, we have $\cG_n\ovee\ssb_{n-1}$ is a  subset of $\cS$. Note that $\cG_n\ovee\ssb_{n-1}$ is finite and $\cS$ is union-closed, so $\ssb_n = \join(\cG_n\ovee\ssb_{n-1}) \in \cS$. Hence
\[ \cG_{n+1} \ovee \ssb_n \subseteq (\cS \ominus \ssb_n) \ovee\ssb_n = \cS\ovee \ssb_n \subseteq \cS, \]
and this completes the inductive step.

Now put $\cF_n = \cG_n\ovee \ssb_{n-1}$ for each $n\in\Nat$. We have just shown that $\cF_n\subseteq \cS$.
Moreover, since $E_n\subseteq\join(\cG_n) \subseteq \Om\setminus \ssb_{n-1}$, $E_n$ is still a witness for~$\cF_n$.

We have $E_{<n} \subseteq \bigsetcup_{j=1}^{n-1} \join(\cG_j) = \ssb_{n-1}\subseteq \meet(\cF_n)$. Finally, if $k > n$, then
\[ E_k \subseteq \join(\cG_k) \subseteq \Om \setminus \ssb_{k-1} \subseteq \Om \setminus \ssb_n \]
so that $E_k \setcap\join(\cF_n) = E_k \setcap \ssb_n = \emptyset$.
\end{proof}
\begin{prop}[Reformulating $\TMIN$]\label{p:reformulate TMIN}
Let $\cS\subseteq \cP(\Om)$ be a union-closed set system, and let $(E_n)_{n\geq 1}$ be a sequence in $\FIN(\Om)$, with  $|E_n|\geq 2$ for all $n$. The following are equivalent:
\begin{romnum}
\item\label{li:TMIN-ish}
there exist incompressible subsets $\cF_n\subseteq \cS$ for all $n\in \Nat$, which satisfy the conditions (WIT) and (MIN) with respect to $(E_n)_{n\ge 1}$;
\item\label{li:TMIN-ind-new}
there exist incompressible subsets $\cH_n\subseteq \cP(\Om)$, for $n\geq 1$, with $E_n$ a witness for $\cH_n$, such that $\cH_1\subseteq \cS$ and
\[
\cH_{n+1}\subseteq \cS^{-E_{<n+1}} \owedge \meet(\cH_n) \qquad \text{ for all } n\in \Nat.
\]
\end{romnum}
\end{prop}
\begin{proof}
\ref{li:TMIN-ish}$\implies$\ref{li:TMIN-ind-new}:
First, recall that (WIT) and (MIN) together imply that $E_j\setcap E_k=\emptyset$ for all $j<k$.
Put $\cH_1=\cF_1\subseteq\cS$ and $\cH_n= \cF_n \owedge \bigcap_{j=1}^{n-1}\meet(\cF_j)$ for all $n\geq 2$.

Since $E_n$ is a witness for $\cF_n$ and is a subset of $\meet(\cF_j)$ for all $j<n$ by condition~(MIN), $E_n$ remains a witness for~$\cH_n$.
Moreover, for every $n\geq 1$ we have $\meet(\cH_n)=\bigcap_{j=1}^n\meet(\cF_j)$. Hence
\[
  \cH_{n+1} =\cF_{n+1} \owedge \meet(\cH_n)
  \subseteq \cS^{-E_{<n+1}} \owedge \meet(\cH_n),
\]
as required.

\medskip
\ref{li:TMIN-ind-new}$\implies$\ref{li:TMIN-ish}:
First note that the conditions on $\cH_n$  imply that
\begin{equation}
\label{eq:easy nest}
\meet(\cH_j)\subseteq \join(\cH_j)  \subseteq \join(\cS^{-E_{<j}} \owedge \meet(\cH_{j-1}))\subseteq \meet(\cH_{j-1})  \qquad\text{for all $j\geq 2$.}
\tag{$\dagger$}
\end{equation}
In particular, if $n\in\Nat$ and $k>n$ then
\begin{equation}
\label{eq:tail nest}
E_k \subseteq \join(\cH_k) \subseteq \meet(\cH_{k-1}) \subseteq \dots \subseteq \meet(\cH_n).
\tag{$\ddagger$}
\end{equation}

Take $\cF_1=\cH_1\subseteq \cS$, noting that $E_1$ is a witness for $\cF_1$. For each $n\geq 2$, since $\cH_n\subseteq \cS^{-E_{<n}}\owedge \meet(\cH_{n-1})$, Lemma~\ref{l:lift inc} ensures any lifting of $\cH_n$ to some $\cF_n\subseteq \cS^{-E_{<n}}$ still has $E_n$ as a witness. This ensures that 
$E_{<n}\cap \join(\cF_n)=\emptyset$  for all $n \geq 2$; while for each $k> n$, the inclusions in \eqref{eq:tail nest} yield $E_k \subseteq \meet(\cH_n)\subseteq\meet(\cF_n)$.
\end{proof}

\begin{prop}[Reformulating $\TORT$]\label{p:reformulate TORT}
Let $\cS\subseteq \cP(\Om)$ be a union-closed set system, and let $(E_n)_{n\geq 1}$ be a sequence in $\FIN(\Om)$, with  $|E_n|\geq 2$ for all $n$. The following are equivalent:
\begin{romnum}
\item\label{li:TORT-ish}
there exist incompressible subsets $\cF_n\subseteq \cS$ for all $n\in \Nat$, which satisfy the conditions (WIT) and (ORT)  with respect to $(E_n)_{n\ge 1}$;
\item\label{li:TORT-ind-new}
there exist incompressible subsets $\cL_n\subseteq \cP(\Om)$, for $n\geq 1$, with $E_n$ a witness for $\cL_n$, such that $\cL_1\subseteq \cS$ and
\[
\cL_{n+1}\subseteq \cS_{E_{<n+1}} \owedge \meet(\cL_n) \qquad\text{for all } n\in\Nat.
\]
\end{romnum}
\end{prop}
\begin{proof}
\ref{li:TORT-ish}$\implies$\ref{li:TORT-ind-new}:
First, recall that (WIT) and (ORT) together imply that $E_j\setcap E_k=\emptyset$ for all $j<k$.
Put $\cL_1=\cF_1\subseteq\cS$ and $\cL_n= \cF_n \owedge \bigcap_{j=1}^{n-1}\meet(\cF_j)$ for all $n\geq 2$.

Since $E_n$ is a witness for $\cF_n$ and is a subset of $\meet(\cF_j)$ for all $j<n$ by condition~(ORT), $E_n$ remains a witness for~$\cL_n$.
Moreover, for every $n\geq 1$ we have $\meet(\cL_n)=\bigcap_{j=1}^n\meet(\cF_j)$. Hence
\[
  \cL_{n+1} =\cF_{n+1} \owedge \meet(\cL_n)
  \subseteq \cS_{-E_{<n+1}} \owedge \meet(\cL_n),
\]
as required.

\medskip
\ref{li:TORT-ind-new}$\implies$\ref{li:TORT-ish}:
First note that the conditions on $\cL_n$ imply that
\begin{equation}
\label{eq:easy nest TORT}
\meet(\cL_j)\subseteq \join(\cL_j) \subseteq \join(\cS_{E_{<j}} \owedge \meet(\cL_{j-1})) \subseteq \meet(\cL_{j-1})  \qquad\text{for all $j\geq 2$.}
\tag{$\dagger$}
\end{equation}
In particular, if $n\in\Nat$ and $k> n$ then
\begin{equation}
\label{eq:tail nest TORT}
E_k \subseteq \join(\cL_k) \subseteq \meet(\cL_{k-1}) \subseteq \dots \subseteq \meet(\cL_n).
\tag{$\ddagger$}
\end{equation}
Take $\cF_1=\cL_1\subseteq \cS$, noting that $E_1$ is a witness for $\cF_1$. For each $n\geq 2$, Lemma~\ref{l:lift inc} ensures that any lifting of $\cL_{n} \subseteq \cS_{E_{<n}}\owedge \meet(\cL_{n-1})$ to some $\cF_n\subseteq \cS_{E_{<n}}$ still has $E_n$ as a witness. This ensures that $E_{<n}\subseteq \meet(\cF_n)$ for all $n \geq 2$; while for each $k> n$, the inclusions in \eqref{eq:tail nest TORT} yield $E_k \subseteq \meet(\cL_n)\subseteq\meet(\cF_n)$.
\end{proof}

\end{subsection}
\end{section}
\begin{section}{Proof of the main theorem}\label{section:pf}
The following definitions are motivated by the strategy outlined in Section \ref{ss:overview of proof}.
From this point onwards, for a set system $\cS\subseteq \cP(\Om)$, we abbreviate ``union-closed with infinite breadth'' to \ucib.

\begin{dfn}[Obstacles]\label{def:JWblock}
Let $\cT$ be \ucib.
\begin{itemize}
\item
We say $\cT$ is \dt{$J$-blocked} if $\exists\ k\geq 2$ so that $\forall\ \cF\in\Inc_k(\cT)$, we have $\br(\cT\ominus\join(\cF))<\infty$.
\item
We say $\cT$ is \dt{$W$-blocked} if $\exists\ r\geq 2$ so that $\forall\ \cF\in\Inc_r(\cT)$ and any witness $F$ for $\cF$, we have $\br(\cT^{-F})<\infty$.
\end{itemize}
\end{dfn}

The idea is that if a set system is $J$-blocked then removing the join of any sufficiently large incompressible set will lead to a dead-end, while if it is $W$-blocked then excluding any sufficiently large witness will also lead to a dead-end.

\begin{lem}[Being blocked passes to subprojections]
\label{l:blocked downwards}
Let $\cT_1,\cT_2$ be \ucib\ set systems on $\Om$ with $\cT_2\subproj \cT_1$.
If $\cT_1$ is $J$-blocked, then so is $\cT_2$.
If $\cT_1$ is $W$-blocked, then so is $\cT_2$.
\end{lem}

\begin{lem}[Step towards the $\TMIN$ case]
\label{l:step for TMIN}
Let $\cT$ be \ucib. If $\cT$ is $J$-blocked but not $W$-blocked, then
\begin{quote}
$\forall\ n\geq 2$ $\exists\ \cH \in \Inc_n(\cT)$ with a witness $H$, such that $\br(\cT^{-H}\owedge\meet(\cH))=\infty$.
\end{quote}
\end{lem}
\begin{lem}[Step towards the $\TORT$ case]
\label{l:step for TORT}
Let $\cT$ be \ucib. If $\cT$ is $J$-blocked and $W$-blocked, then
\begin{quote}
$\forall\ n\geq 2$ $\exists\ \cL \in \Inc_n(\cT)$ with a witness $L$, such that $\br(\cT_L\owedge\meet(\cL))=\infty$.
\end{quote}
\end{lem}
The proofs of these lemmas are deferred to the next section; they are not difficult, but the details would interrupt the flow of the argument. Taking these lemmas on trust for now, we use the concepts of $J$-blocked and $W$-blocked to define three natural conditions, each of which will ensure that a subprojection of the form $\TMAX$, $\TMIN$ or $\TORT$ may be constructed by ``following one's nose''.

\begin{dfn}[Conditions ensuring generation of $\TMAX$, $\TMIN$ or $\TORT$]
\label{d:automatic}
Let $\cS$ be \ucib.
\begin{romnum}
\item
\label{li:tmax-auto}
We say $\cS$ is \dt{$\TMAX$-automatic} if every \ucib\ $\cT\subproj\cS$ is not $J$-blocked.

\item
\label{li:tmin-auto}
We say $\cS$ is \dt{$\TMIN$-automatic} if every \ucib\ $\cT\subproj\cS$ is $J$-blocked but not $W$-blocked.

\item
\label{li:tort-auto}
We say $\cS$ is \dt{$\TORT$-automatic} if it is $J$-blocked and $W$-blocked.
(By Lemma \ref{l:blocked downwards}, this is equivalent to requiring that every \ucib\ $\cT\subproj\cS$ is $J$-blocked and $W$-blocked.)
\end{romnum}
\end{dfn}

%
\begin{prop}[Ensuring a subprojection which is $\cT_*$-automatic]
\label{p:auto-inevitable}
Let $\cS$ be \ucib. Suppose $\cS$ is not $\TMAX$-automatic. Then $\cS$ has a \ucib\ subprojection which is $\TMIN$-automatic or one which is $\TORT$-automatic.
\end{prop}

\begin{proof}
Since $\cS$ is not $\TMAX$-automatic, there exists some \ucib\ $\cS_0 \subproj \cS$ which is $J$-blocked. Note that by Lemma \ref{l:blocked downwards} every subprojection of $\cS_0$ is also $J$-blocked.
Now, there are two cases to consider.

\begin{itemize}
\item \textit{Case 1: every \ucib\ $\cT\subproj\cS_0$ is not $W$-blocked.}
By the remark above, every such $\cT$ is also $J$-blocked, and so in this case $\cS_0$ is $\TMIN$-automatic with $\cS_0\subproj\cS$.
\item \textit{Case 2: there is some \ucib\ $\cS_1\subproj \cS_0$ which is $W$-blocked.}
By the remark above, $\cS_1$ is also $J$-blocked. So in this case, $\cS_1$ is $\TORT$-automatic with $\cS_1\subproj\cS$ (by transitivity).
\end{itemize}
This completes the proof.
\end{proof}
\begin{prop}[Justifying Definition \ref{d:automatic}]
\label{p:auto-generate}
Let $\cS$ be \ucib, and let $\cT_*$ denote either $\TMAX$, $\TMIN$ or $\TORT$. If $\cS$ is $\cT_*$-automatic, then there exists
$\cE$ as in Definition~\ref{d:3kings}
such that $\cT_*(\cE)\subproj\cS$.
\end{prop}

\begin{proof}
For this proof, we combine Lemma \ref{l:step for TMIN} and \ref{l:step for TORT} with the ``inductive reformulations'' of $\TMAX$, $\TMIN$ and $\TORT$.
We go through each case in turn.

\paragraph{The $\TMAX$ case.}
Suppose $\cS$ is $\TMAX$-automatic.
Put $\cT_0=\cS$; then we may inductively construct sequences $(E_n)_{n\geq 1}$ and $(\cG_n)_{n\geq 1}$, such that the following properties hold.
\begin{itemize}
\item $E_n\subseteq\Om$ and $|E_n|=n+1$; 
\item $E_n$ is a witness for some incompressible $\cG_n\subseteq \cT_{n-1}$;
\item $\cT_n\defeq \cT_{n-1} \ominus \join(\cG_n)$  has infinite breadth.
\end{itemize}
At each stage, $\cT_n \subproj\cT_{n-1}\subproj \dots \subproj \cT_0 =\cS$; so the existence of suitable $\cG_n$ and $E_n$ is guaranteed by the definition of $\TMAX$-automaticity.

The recursive construction of $\cT_n$ yields $\cT_n = \cS \ominus\left( \bigsetcup_{j=1}^{n-1} \join(\cG_j)\right)$. 
Hence, by the direction \ref{li:TMAX-ind-new}$\implies$\ref{li:TMAX-ish} in Proposition \ref{p:reformulate TMAX},  together with Lemma \ref{l:TMAX-orig}, $\cS\owedge \left(\bigsetcup_{n\geq 1} E_n\right)$ contains the configuration $\TMAX(\cE)$.

\paragraph{The $\TMIN$ case.}
Suppose $\cS$ is $\TMIN$-automatic. Put $\cT_0=\cS$; then we may inductively construct sequences $(E_n)_{n\geq 1}$ and $(\cH_n)_{n\geq 1}$, such that the following properties hold.
\begin{itemize}
\item $E_n\subseteq\Om$ and $|E_n|=n+1$; 
\item $E_n$ is a witness for some incompressible $\cH_n\subseteq \cT_{n-1}$;
\item $\cT_n\defeq (\cT_{n-1})^{-E_n} \owedge \meet(\cH_n)$  has infinite breadth.
\end{itemize}
At each stage, $\cT_n \subproj\cT_{n-1}\subproj \dots \subproj \cT_0 =\cS$; so the existence of such $\cH_n$ and $E_n$ is guaranteed by Lemma~\ref{l:step for TMIN}.

To prove that $\cS\owedge \left(\bigsetcup_{n\geq 1} E_n\right)$ contains $\TMIN(\cE)$, it suffices to show that $(\cH_n)$ and $(E_n)$ satisfy the conditions in Proposition \ref{p:reformulate TMIN}\ref{li:TMIN-ind-new}, for then we can apply that proposition together with Lemma \ref{l:TMIN-orig}. This will follow from the properties above and the following claim.

\paragraph{Claim:} $\cT_j = \cS^{-E_{<j+1}} \owedge\meet(\cH_j)$ for all $j\in\Nat$.

The claim holds when $j=1$, since
\[ \cT_1 = (\cT_0)^{-E_1} \owedge \meet(\cH_1) = \cS^{-E_{<2}} \owedge \meet(\cH_1). \]
Assume it holds for some $j\in\Nat$; then
\[ \begin{aligned}
\cT_{j+1}
& = (\cT_j)^{-E_{j+1}} \owedge \meet(\cH_{j+1}) & \text{by definition}\\
& = \left( \cS^{-E_{<j+1}} \owedge\meet(\cH_j) \right)^{-E_{j+1}} \owedge  \meet(\cH_{j+1}) & \text{by the inductive hypothesis} \\
& = \cS^{-(E_{<j+1}\setcup E_{j+1})} \owedge \left( \meet(\cH_j) \setcap \meet(\cH_{j+1})\right)  &  \text{by Lemma \ref{l:swap-proj}\ref{li:swap exclude}} \\
& = \cS^{-E_{<j+2}} \owedge  \meet(\cH_{j+1}),
\end{aligned} \]
as required.
(The third and fourth steps are justified by the inclusions in \eqref{eq:easy nest}. Specifically, we can apply Lemma \ref{l:swap-proj} because
the inductive hypothesis ensures that $E_{j+1}\subseteq\join(\cH_{j+1})\subseteq\meet(\cH_j)$; and we also have $\meet(\cH_{j+1})\subseteq\join(\cH_{j+1})\subseteq\meet(\cH_j)$.)
By induction, the claim holds for all $j\in\Nat$.

\paragraph{The $\TORT$ case.}
Suppose $\cS$ is $\TORT$-automatic. Put $\cT_0=\cS$; then we may inductively construct sequences $(E_n)_{n\geq 1}$ and $(\cL_n)_{n\geq 1}$, such that the following properties hold.
\begin{itemize}
\item $E_n\subseteq\Om$ and $|E_n|=n+1$; 
\item $E_n$ is a witness for some incompressible $\cL_n\subseteq \cT_{n-1}$;
\item $\cT_n\defeq (\cT_{n-1})_{E_n} \owedge \meet(\cL_n)$  has infinite breadth.
\end{itemize}
At each stage, $\cT_n \subproj\cT_{n-1}\subproj \dots \subproj \cT_0 =\cS$; so the existence of suitable $\cL_n$ and $E_n$ is guaranteed by Lemma~\ref{l:step for TORT}.

To prove that $\cS\owedge \left(\bigsetcup_{n\geq 1} E_n\right)$ contains $\TORT(\cE)$, it suffices to show that $(\cL_n)$ and $(E_n)$ satisfy the conditions in Proposition \ref{p:reformulate TORT}\ref{li:TORT-ind-new}; for then we can apply that proposition together with Lemma \ref{l:TORT-orig}. This will follow from the properties above and the following claim.

\paragraph{Claim:} $\cT_j = \cS_{E_{<j+1}} \owedge\meet(\cL_j)$ for all $j\in\Nat$.

The claim holds when $j=1$, since
\[ \cT_1 = (\cT_0)_{E_1} \owedge \meet(\cL_1) = \cS_{E_{<2}} \owedge \meet(\cL_1). \]
Assume it holds for some $j\in\Nat$; then
\[ \begin{aligned}
\cT_{j+1}
& = (\cT_j)_{E_{j+1}} \owedge \meet(\cL_{j+1}) & \text{by definition}\\
& = \left( \cS_{E_{<j+1}} \owedge\meet(\cL_j) \right)_{E_{j+1}} \owedge  \meet(\cL_{j+1}) & \text{by the inductive hypothesis} \\
& = \cS_{(E_{<j+1}\setcup E_{j+1})} \owedge \left( \meet(\cL_j) \setcap \meet(\cL_{j+1})\right)  &  \text{by Lemma \ref{l:swap-proj}\ref{li:swap include}} \\
& = \cS_{E_{<j+2}} \owedge  \meet(\cL_{j+1}),
\end{aligned} \]
as required.
(The third and fourth steps are justified by the inclusions in \eqref{eq:easy nest TORT}. Specifically, we can apply Lemma \ref{l:swap-proj} because the inductive hypothesis ensures that
 $E_{j+1}\subseteq\join(\cL_{j+1})\subseteq\meet(\cL_j)$; and we also have $\meet(\cL_{j+1})\subseteq\join(\cL_{j+1})\subseteq\meet(\cL_j)$.)
By induction, the claim holds for all $j\in\Nat$.

\medskip
This completes the proof of the proposition.
\end{proof}

\paragraph{Putting everything together.}
Combining Proposition \ref{p:auto-generate} and Proposition \ref{p:auto-inevitable}, and using transitivity of $\subproj$ again, we deduce that a \ucib\ set system $\cS\subseteq\cP(\Om)$ always has a subprojection of the form $\TMAX(\cE)$, $\TMIN(\cE)$ or $\TORT(\cE)$.

Therefore, to complete the proof of Theorem \ref{t:headline result}, it remains only to prove Lemmas \ref{l:blocked downwards}, \ref{l:step for TMIN} and \ref{l:step for TORT}. This will be dealt with in the next section.
\end{section}

\begin{section}{Proofs of the technical lemmas}\label{section:lemmas}
\begin{proof}[Proof of Lemma \ref{l:blocked downwards} (being blocked passes to subprojections)]
\strut\newline
 First, let $\cU\subseteq\cT_1$ and $X\subseteq\Om$ be such that $\cT_2=\cU\owedge X$.
Suppose $\cT_1$ is $J$-blocked, and let $k\geq 2$ be as in Definition \ref{def:JWblock}.
Let $\cF\in\Inc_k(\cT_2)$. By Lemma~\ref{l:lift inc},  any lifting $\cG\subseteq \cU$ of $\cF$ satisfies $\cG\in\Inc_k(\cU)$. Also, since $\cG\owedge X = \cF$, we have $\join(\cG)\setcap X =\join(\cF)$. It follows that
\[ \begin{aligned}
\cT_2 \ominus \join(\cF)
 & = (\cU\owedge X)\ominus \join(\cF) \\
& = (\cU \owedge X) \ominus \join(\cG) \\
& = (\cU \ominus \join(\cG)) \owedge X & \subproj \cT_1\ominus \join(\cG).
\end{aligned} \]
Since taking subprojections cannot increase breadth
(Proposition \ref{p:breadth of subproj}),
and since $\cT_1$ is $J$-blocked ``at depth $k$'', it follows that
$\br(\cT_2\ominus\join(\cF))<\infty$.
Since this holds for all $\cF\in\Inc_k(\cT_2)$, $\cT_2$ is $J$-blocked.

Similarly, assume $\cT_1$ is $W$-blocked, and let $r\geq 2$ be as in Definition \ref{def:JWblock}. 
Let $\cF\in\Inc_r(\cT_2)$ and let $F$ be any witness for $\cF$. By Lemma~\ref{l:lift inc}
any lifting of $\cF$ to some $\cG\subseteq \cU$  admits $F$ as a witness as well, and $\cG\in\Inc_r(\cU)$. Observe that
\[ \cT_2^{-F} = (\cU\owedge X)^{-F} = \cU^{-F} \owedge X \subseteq \cT_1^{-F}\owedge X,
\]
where the second equality follows from
Lemma~\ref{l:swap-proj}\ref{li:swap exclude} and the fact that $F\subseteq \join(\cT_2)\subseteq X$.
In~particular, $\cT_2^{-F}\subproj \cT_1^{-F}$.
Since taking subprojections cannot increase breadth, and since $\cT_1$ is $W$-blocked ``at depth $r$'', we have $\br(\cT_2^{-F}) <\infty$. Since this holds for any witness of any $\cF\in\Inc_r(\cT_2)$, $\cT_2$ is $W$-blocked.
\end{proof}

\begin{proof}[Proof of Lemma~\ref{l:step for TMIN} (step towards the $\TMIN$ case)]
\strut\newline
We start by fixing some $k\geq 2$ such that
$\br(\cT\ominus\join(\cF))<\infty$ for all $\cF\in\Inc_k(\cT)$. (Such a $k$ exists since $\cT$ is $J$-blocked.)

Now, let $n\geq 2$. Since $\cT$ is not $W$-blocked, there exists $\cG\in\Inc_{kn}(\cT)$ with a witness $G$ such that $\br(\cT^{-G})=\infty$.
Enumerate $\cG$ as $\{ \ssx_{ij} \colon 1\leq i\leq k, 1\leq j\leq n\}$ and enumerate $G$ as $\{\gamma_{ij}\colon 1\leq i\leq k, 1\leq j\leq n\}$, in such a way that $\gamma_{ij}\in \ssx_{ij}$ for all $i$,$j$.  

For each $j$, let $\cF_j=\{\ssx_{ij}\colon 1\leq i\leq k\}$ and let $\ssy_j=\join(\cF_j)\in\cT$.
Then put $H = \{ \gamma_{1,j}\colon 1\leq j\leq n\}$ and let $\cH=\{ \ssy_1,\dots,\ssy_n\}\subseteq\cT$. By the coarse-graining lemma (Lemma~\ref{l:coarse-grain}), $\cH\in\Inc_n(\cT)$ with $H$ as a witness.

Since $H\subseteq G$, we have $\cT^{-H} \supseteq \cT^{-G}$, and so $\br(\cT^{-H}) \geq \br(\cT^{-G}) = \infty$. On the other hand, since $\cF_j\in\Inc_k(\cT)$ and $\join(\cF_j)=\ssy_j$, we have $\br(\cT\ominus\ssy_j)<\infty$, and so $\br(\cT^{-H}\ominus\ssy_j)<\infty$ for $j=1,\dots, n$. Applying Lemma \ref{l:swissarmyknife}\ref{li:breadth of partition} to $\cT^{-H}$ and 
\[
	\Omega=\meet(\cH)\cup\bigcup_{j=1}^n(\Omega\setminus \ssy_j)
\] 
we conclude that   $\br(\cT^{-H}\owedge \meet(\cH))=\infty$.
\end{proof}

\begin{proof}[Proof of Lemma~\ref{l:step for TORT} (step towards the $\TORT$-case)]
  \strut\newline
  Fix $k,r\geq 2$ such that:
\begin{equation}\label{eq:J-blocked}
\forall\ \cF\in\Inc_k(\cT) \quad \br(\cT\ominus\join(\cF))<\infty ;
\tag{$*$}
\end{equation}
\begin{equation}\label{eq:W-blocked}
\forall\ \cF\in\Inc_r(\cT) \text{ and any witness $F$ for $\cF$}, \quad \br(\cT^{-F})<\infty\;.
\tag{$**$}
\end{equation}
(Such $k$ and $r$ exist by the assumption that $\cT$ is $J$-blocked and $W$-blocked.)

Now, let $n\geq 2$, and pick any $\cG\in\Inc_{k(n+r-1)}(\cT)$. Enumerate $\cG$ as $\{\ssz_{ij} \colon 1\leq i\leq k, 1\leq j \leq n+r-1\}$ and pick a corresponding witness $G=\{\gamma_{ij} \colon 1\leq i\leq k, 1\leq j \leq n+r-1\}$ with the property that $\gamma_{ij} \in\ssz_{ij}$ for all $i,j$.

Let $V\defeq \{ \gamma_{1,j} \colon 1\leq j \leq r+n-1\}$. For any $\ssw\subseteq\Om$ note that either $|\ssw\setcap V | \geq n$ or $|V\setminus \ssw | \geq r$ (pigeon-hole principle); hence, either $\ssw$ contains some $n$-element subset of $V$, or it is disjoint from some $r$-element subset of $V$. Hence, for any set system $\cJ \subseteq\cP(\Om)$,
\begin{equation}\label{eq:include or exclude}
\cJ = \left( \bigsscup_{L\subseteq V\colon |L|=n} \cJ_L \right) \sscup \left( \bigsscup_{F\subseteq V \colon |F| = r} \cJ^{-F} \right)
\tag{$***$}
\end{equation}

We apply this to $\cT$.
By the condition \eqref{eq:W-blocked}, we have $\br(\cT^{-F})<\infty$ for all $r$-element subsets $F\subseteq V$.
(See Lemma \ref{l:easy witness}\ref{li:subset of witness}.)
On the other hand, $\br(\cT)=\infty$ by assumption. Therefore,
by Lemma \ref{l:swissarmyknife}\ref{li:breadth of merged}
and the identity \eqref{eq:include or exclude}, we must have $\br(\cT_L)=\infty$ for some $n$-element subset $L\subseteq V$. Pick such an $L$; by relabelling if necessary, we may assume that
\[ L = \{ \gamma_{1,j} \colon 1 \leq j \leq n \}. \]

Now let $\ssx_j = \bigsetcup_{i=1}^k \ssz_{ij}$, for $1\leq j \leq n$.
Since $\cG$ is incompressible so are all its subsets (this follows from Lemma \ref{l:easy witness}\ref{li:subset of witness}). Hence
for each $j$, $\{\ssz_{ij} \colon 1\leq i \leq k\} \in \Inc_k(\cT)$;
and so, by the condition \eqref{eq:J-blocked}, $\br(\cT\ominus \ssx_j)<\infty$.
Moreover, if we put
\[ \cL\defeq \{ \ssx_j \colon 1 \leq j \leq n\} \]
then $L$ is a witness for $\cL$, by the coarse-graining lemma (Lemma~\ref{l:coarse-grain}).

Finally,
$\br(\cT_L)=\infty$ while $\br(\cT_L\ominus \ssx)<\infty$ for each $\ssx\in\cL$.
Applying Lemma \ref{l:swissarmyknife}\ref{li:breadth of partition} to $\cT_L$ and 
\[
\Omega=\meet(\cL)\cup\bigcup_{j=1}^n(\Omega\setminus \ssx_j)
\]
we conclude that $\br(\cT_L\owedge \meet(\cL))=\infty$.
\end{proof}

\end{section}

\section*{Acknowledgements}
This work grew out of conversations between the authors while attending the conference ``Banach Algebras and Applications'', held in Gothenburg, Sweden, July--August 2013, and was further developed while the authors were attending the thematic program ``Abstract Harmonic Analysis, Banach and Operator Algebras''  at the Fields Institute, Canada, during March--April 2014. The authors thank the organizers of these meetings for invitations to attend and for pleasant environments to discuss research.

The first author acknowledges the financial support of the Faculty of Science and Technology at Lancaster University, in the form of a travel grant to attend the latter meeting. The second author acknowledges financial support from the University of Delaware Research Foundation. The third author acknowledges the financial supports of a Fast Start Marsden Grant and of Victoria University of Wellington to attend these meetings.

The final write-up of this paper was completed during a visit of the first author in March 2020 to the University of Delaware. He thanks the Department of Mathematics and Statistics at Lancaster University for financial support through their Visitor Fund. He also thanks the Department of Mathematical Sciences at the University of Delaware for their hospitality.
The second author also acknowledges support from National Science Foundation grant DMS-1902301 during the preparation of this article.

The authors thank two anonymous referees for their useful advice and corrections, which have led to significant improvements of the presentation of this paper.

\newcommand{\etalchar}[1]{$^{#1}$}

\vfill

\newcommand{\address}[1]{{\small\sc#1.}}
\newcommand{\email}[1]{\texttt{#1}}

\noindent
\address{Yemon Choi,
Department of Mathematics and Statistics,
Lancaster University,
Lancaster LA1 4YF, United Kingdom} 

\email{y.choi1@lancaster.ac.uk}

\noindent
\address{Mahya Ghandehari,
Department of Mathematical Sciences,
University of Delaware,
Newark, Delaware 19716, United States of America}

\email{mahya@udel.edu}

\noindent
\address{Hung Le Pham,
School of Mathematics and Statistics,
Victoria University of Well\-ing\-ton,
Wellington 6140, New Zealand}

\email{hung.pham@vuw.ac.nz}
\end{document}